\documentclass[11pt]{amsart}
\usepackage{amsmath,amssymb,amsthm}
\usepackage{amsfonts}
\usepackage[mathscr]{eucal}
\usepackage{amsmath,amssymb,amsthm}
\usepackage{amsfonts}
\usepackage[mathscr]{eucal}
\newtheorem{theorem}{Theorem}[section]

\newtheorem{corollary}[theorem]{Corollary}

\theoremstyle{definition}

\theoremstyle{remark}
\newtheorem{remark}[theorem]{Remark}

\numberwithin{equation}{section}

\newcommand{\ep}{\varepsilon}

\def\R{{\mathbb R}}
\def\N{{\mathbb N}}

\def\a{\alpha}

\def\ep{\varepsilon}

\def\E{{\mathbb E}}
\def\P{{\mathbb P}}

\def\N{{\mathbb N}}

\begin{document}

\title{\bf A Strong Law of Large Numbers with Applications to
Self-Similar Stable Processes}

\author{Erkan Nane}
\address{Erkan Nane, Department of Mathematics and Statistics,
Auburn University, Auburn, AL 36849}
\email{nane@auburn.edu}
\urladdr{http://www.auburn.edu/$\sim$ezn0001}

\author{Yimin Xiao}
\address{Yimin Xiao, Department Statistics and Probability,
Michigan State University, East Lansing, MI 48823}
\email{xiao@stt.msu.edu}
\urladdr{http://www.stt.msu.edu/$\sim$xiaoyimi}
\thanks{Research of Y. Xiao was partially supported by
NSF grant DMS-0706728.}

\author{Aklilu Zeleke}
\address{Aklilu Zeleke, Department of  Statistics and Probability,
Michigan State University, East Lansing, MI 48823}
\email{zeleke@stt.msu.edu}

\begin{abstract}
Let $p \in (0, \infty)$ be a constant and let $\{\xi_n\} \subset L^p(\Omega,
{\mathcal F}, \P)$ be a sequence of random variables. For any integers
$m, n \ge 0$, denote $S_{m, n} = \sum_{k=m}^{m + n} \xi_k$. It is proved
that, if there exist a nondecreasing function $\varphi: \R_+\to \R_+$
(which satisfies a mild regularity condition) and an
appropriately chosen integer $a\ge 2$ such that
$$
\sum_{n=0}^\infty \sup_{k \ge 0} \E\bigg|\frac{S_{k, a^n}} {\varphi(a^n)}
\bigg|^p < \infty,$$
Then
$$
\lim_{n \to \infty} \frac{S_{0, n}} {\varphi(n)} = 0\qquad \hbox{a.s.}
$$
This extends Theorem 1 in Levental, Chobanyan and Salehi \cite{chobanyan-l-s}
and can be applied conveniently to a wide class of self-similar processes with
stationary increments including stable processes.
\end{abstract}


\maketitle


\section{Introduction}

The study on strong law of large numbers has a long history and there is a
vast body of references on this topic. This note is motivated by our interest in
studying asymptotic properties of stochastic processes with heavy-tailed distributions.
Typical examples of such processes are linear fractional stable motion and harmonizable
fractional stable motion. See Samorodnitsky and Taqqu \cite{ST94} and Embrechts
and Maejima \cite{EmbrechtsMaejima}.

Let $p \in (0, \infty)$ be a constant and let $\{\xi_n\} \subset L^p(\Omega,
{\mathcal F}, \P)$ be a sequence of random variables. For any integers
$m, n \ge 0$, denote
$$
S_{m, n} = \sum_{k=m}^{m + n} \xi_k, \qquad M_{m, n} = \max_{k \le n} |S_{m, k}|.
$$
For any nondecreasing function $\varphi: \R_+\to \R_+$ such that $\varphi(x) \uparrow \infty$
as $x \to \infty$, we say that $\{\xi_n\}$ satisfies the SLLN with respect to $\varphi$
(or $\varphi$-SLLN) if
\[
\lim_{n \to \infty} \frac{S_{0, n}} {\varphi(n)} = 0\qquad \hbox{a.s.}
\]

The following is the main result of this paper, which is an extension of Theorem 1 in
Levental, Chobanyan and Salehi \cite{chobanyan-l-s} who considered the case of $p > 1$
and $\varphi(n) = n$. Our result improves their theorem in two aspects: (i) it gives sharper
norming constants (see Remark 1.2 below) and (ii) it can be applied conveniently to a wide
class of self-similar processes with stationary increments such as self-similar
$\alpha$-stable processes for all $\alpha \in (0, 2]$.

\begin{theorem}\label{Th:SLLN}
Let $p>0$ be a constant and let $\varphi: \R_+\to \R_+$ be a nondecreasing function
such that $\varphi(x)\uparrow \infty$ as $x \to \infty$ and $C_1 \le \varphi(2x) /\varphi(x)
\leq C_2$ for all $x\in \R_+$ for some constants $C_2 \ge C_1 > 1$. Assume  $a\geq 2$
is an integer that satisfies  $C_1^{ p \lfloor\log_2 a\rfloor} \ge \max\{2, 2^p\}$. If we have
\begin{equation}\label{Eq:MomCon1}
\sum_{n=0}^\infty \sup_{k \ge 0} \E\bigg|\frac{S_{k, a^n}} {\varphi(a^n)}
\bigg|^p < \infty,
\end{equation}
then $\{\xi_n\}$ satisfies the $\varphi$-SLLN.
\end{theorem}

\begin{remark}
If $\varphi(n)=n^q$ ($q > 0$) or $\varphi(n)=n^q (\log n)^\beta$ for $q > 0$
and $\beta \in \R$, then $\varphi$ satisfies the conditions of Theorem \ref{Th:SLLN}. From
Theorem \ref{Th:SLLN} it is easy to see that
if  $\{\xi_n, n \ge 1\}$ is a sequence of i.i.d. random variables with mean 0 and variance 1,
then for any $\ep > 0$
\begin{equation}\label{Eq:c-1}
\lim_{n \to \infty} \frac{S_{0,n}} {\sqrt{n \log n (\log \log n)^{1 + \ep}} }= 0\qquad
\hbox{ a.s.}
\end{equation}
Information on higher moments of $S_{k,n}$ leads to improvement on the power of $\log n$.
For example, if $\{\xi_n, n \ge 1\}$ are i.i.d. standard normal random variables, then
for any $\ep > 0$
\begin{equation}\label{Eq:c-2}
\lim_{n \to \infty} \frac{S_{0,n}} {\sqrt{n}\, (\log n)^\ep }= 0\qquad \hbox{ a.s.}
\end{equation}
Even though (\ref{Eq:c-1}) and (\ref{Eq:c-2}) are less precise than the law of the
iterated logarithm, the advantage of this method is that it is still applicable
when the random variables $\{\xi_n, n \ge 1\}$ are dependent or non-Gaussian.
\end{remark}

\noindent{\it Proof of Theorem \ref{Th:SLLN}}
The proof of Theorem \ref{Th:SLLN} is a modification of the proof
of Theorem 1 in Levental, Chobanyan and Salehi \cite{chobanyan-l-s}.

Let $\N_0=\N\cup \{0\}$. We first consider $p>1$.
For any $k\in \N_0, n\in\N_0$
\begin{equation}\label{Eq:Mkp}
M_{k, a^{n+1}}\leq \max\{M_{k, a^{n}}, |S_{k,a^n}|+M_{k+a^n, a^{n}}\}.
\end{equation}
Using the elementary inequality $|x+y|^p\leq 2^{p-1}(|x|^p+|y|^p)$  we get
\begin{eqnarray}
M^p_{k, a^{n+1}}&\leq& \max\Big\{M^p_{k, a^{n}}, 2^{p-1}(|S_{k,a^n}|^p
+M^p_{k+a^n, a^{n}})\Big\} \nonumber\\
&\leq & (2^{p-1}-1)|S_{k,a^n}|^p+M^p_{k, a^{n}}+2^{p-1}M^p_{k+a^n, a^{n}}.
\label{max-1}
\end{eqnarray}
Eq. \eqref{max-1} can be written as
\begin{eqnarray}
M^p_{k, a^{n+1}}-|S_{k,a^{n+1}}|^p &\leq &  M^p_{k, a^{n}}-|S_{k,a^n}|^p
+2^{p-1}\Big(M^p_{k+a^n, a^{n}}-|S_{k+a^n,a^n}|^p\Big)\nonumber\\
& & -|S_{k,a^{n+1}}|^p + 2^{p-1}|S_{k,a^n}|^p+2^{p-1}|S_{k+a^n,a^n}|^p.
\nonumber
\end{eqnarray}
Dividing both sides by $\varphi(a^{n+1})^p$, taking expectations, and
then the supremum over all $k$'s, we get
\begin{eqnarray}
F_{n+1}&\leq& \frac{\varphi(a^{n})^p}{\varphi(a^{n+1})^p}F_n+
\frac{2^{p-1}\varphi(a^{n})^p}{\varphi(a^{n+1})^p}F_n +G_n \nonumber\\
&=& \frac{\varphi(a^{n})^p}{\varphi(a^{n+1})^p}(1+2^{p-1})F_n+G_n, \qquad
 \  n\in \N_0,\label{recursion}
\end{eqnarray}
where
\begin{equation}\label{Eq:Fn}
F_n=\sup_{k\in \N_0}\E\bigg(\frac{M^p_{k, a^{n}}-|S_{k,a^n}|^p}{\varphi(a^{n})^p}\bigg)
\end{equation}
and
\[
\begin{split}
G_n&=\sup_{k\in \N_0}\Bigg\{\frac{2^{p-1}\varphi(a^{n})^p}{\varphi(a^{n+1})^p}
\E\bigg(\bigg|\frac{S_{k+a^n,a^n}}{\varphi(a^{n})}\bigg|^p\bigg)\\
& \qquad \qquad \qquad +
\frac{2^{p-1}\varphi(a^{n})^p}{\varphi(a^{n+1})^p}\E\bigg(\bigg|\frac{S_{k,a^n}}
{\varphi(a^{n})}\bigg|^p\bigg)
- \E\bigg(\bigg|\frac{S_{k,a^{n+1}}}{\varphi(a^{n+1})}\bigg|^p\bigg)\Bigg\}.
\end{split}
\]

Now by the assumptions on $\varphi$ and $a$, we have
\begin{equation}\label{crucial-bound}
\frac{\varphi(a^{n})^p}{\varphi(a^{n+1})^p}(1+2^{p-1})\leq
\frac{1+2^{p-1}} {C_1^{p \lfloor \log_2 a\rfloor}}
:= c<1, \ \qquad n\in \N_0.
\end{equation}
Then $F_{n+1} \le c F_n + G_n$ and it is easy to show by induction in $n$ that
$$
F_{n+1}\leq \sum_{k=0}^n c^{n-k}G_k, \qquad \ n\in \N_0.
$$

By summing up \eqref{recursion} from $n=0$ to $\infty$ we get
\begin{equation}\label{good-upper-bound}
\sum_{n=0}^\infty \E\bigg(\frac{M^p_{a^n, a^{n}}-|S_{a^n,a^n}|^p}{\varphi(a^{n})^p}\bigg)
\leq \frac{1}{1-c}\sum_{n=0}^\infty G_n.
\end{equation}
It follows from (\ref{Eq:MomCon1}) that $\sum_{n=0}^\infty G_n < \infty$.

Next we consider the case $0<p\leq 1$ and establish an inequality similar to
(\ref{good-upper-bound}).
Using (\ref{Eq:Mkp}) and the elementary inequality $|x+y|^p\leq |x|^p+|y|^p$ as in
\cite{chobanyan-l-s} we get
\begin{eqnarray}
M^p_{k, a^{n+1}}&\leq& \max \big\{M^p_{k, a^{n}}, \,
|S_{k,a^n}|^p+M^p_{k+a^n, a^{n}}\big\}\nonumber\\
&\leq & M^p_{k, a^{n}}+M^p_{k+a^n, a^{n}}.\label{max-2}
\end{eqnarray}
It follows that
\begin{eqnarray}
M^p_{k, a^{n+1}}-|S_{k,a^{n+1}}|^p &\leq &  M^p_{k, a^{n}}-
|S_{k,a^n}|^p  +(M^p_{k+a^n, a^{n}}
-|S_{k+a^n,a^n}|^p)\nonumber\\
& & -|S_{k,a^{n+1}}|^p + |S_{k,a^n}|^p+|S_{k+a^n,a^n}|^p.
\nonumber
\end{eqnarray}
Dividing both sides by $\varphi(a^{n+1})^p$, taking expectations, and then
the supremum over all $k$'s, we derive
\begin{eqnarray}
F_{n+1}&\leq& \frac{\varphi(a^{n})^p}{\varphi(a^{n+1})^p}F_n+\frac{\varphi(a^{n})^p}
{\varphi(a^{n+1})^p}F_n +H_n
\nonumber\\
&=& \frac{2\varphi(a^{n})^p}{\varphi(a^{n+1})^p}F_n+H_n, \
\qquad  n\in \N_0,\label{recursion-1}
\end{eqnarray}
where $F_n$ is defined as in (\ref{Eq:Fn}) and
$$
H_n=\sup_{k\in \N_0}\Bigg\{\frac{\varphi(2^{n})^p}{\varphi(a^{n+1})^p}
\E\bigg(\bigg|\frac{S_{k+a^n,a^n}}{\varphi(a^{n})}\bigg|^p\bigg)+
\frac{\varphi(a^{n})^p}{\varphi(a^{n+1})^p}
\E\bigg(\bigg|\frac{S_{k,a^n}}{\varphi(a^{n})}\bigg|^p\bigg)-
\E\bigg(\bigg|\frac{S_{k,a^{n+1}}}{\varphi(a^{n+1})}\bigg|^p\bigg)\Bigg\}.
$$

Again, it follows from the assumptions on $\varphi$ and $a$ that
\begin{equation}\label{crucial-bound-1}
\frac{2\varphi(a^{n})^p}{\varphi(a^{n+1})^p}\leq \frac{2}
{C_1^{p\lfloor \log_2 a\rfloor}}:=c<1, \qquad \ n\in \N_0.
\end{equation}
By using induction in $n$ we derive
$$
F_{n+1}\leq \sum_{k=0}^n c^{n-k}H_k, \quad \ n\in \N_0.
$$
Hence, by summing up \eqref{recursion-1} from $n=0$ to $\infty$ we get
\begin{equation}\label{good-upper-bound-1}
\sum_{n=0}^\infty \E\bigg(\frac{M^p_{a^n, a^{n}}-|S_{a^n,a^n}|^p}{\varphi(a^{n})^p}\bigg)
\leq \frac{1}{1-c}\sum_{n=0}^\infty H_n.
\end{equation}
Eq. (\ref{Eq:MomCon1}) implies that  $\sum_{n=0}^\infty H_n < \infty$.

By combining  \eqref{good-upper-bound} and \eqref{good-upper-bound-1}, we see
that almost surely
\begin{equation}
\frac{M^p_{a^n, a^{n}}-|S_{a^n,a^n}|^p}{\varphi(a^{n})^p}\to 0, \ \ \  \mathrm{as} \ n\to\infty.
\end{equation}
Note that \eqref{Eq:MomCon1} also implies that
\begin{equation}
\frac{|S_{a^n,a^n}|}{\varphi(a^{n})}\to 0, \ \mathrm{as} \ \ \ \ n\to\infty.
\end{equation}
Therefore, we obtain almost surely
\begin{equation}
\frac{M_{a^n, a^{n}}}{\varphi(a^{n})}\to 0, \ \mathrm{as} \ n\to\infty.
\end{equation}
Since
\begin{equation}\label{phi-2}
\inf _n\frac{\varphi(a^{n+1})}{\varphi(a^{n})}\ge C_1^{\lfloor \log_2 a\rfloor} >1
\end{equation}
we have
\begin{equation}
\begin{split}
\frac{M_{a^n, a^{n}}}{\varphi(a^{n})}&=\frac{\varphi(a^{n+1})-\varphi(a^{n})}
{\varphi(a^{n})}\frac{M_{a^n, a^{n}}}{\varphi(a^{n+1})-\varphi(a^{n})}\\
&\geq
(C_1^{\lfloor \log_2 a\rfloor}-1)\frac{M_{a^n, a^{n}}}{\varphi(a^{n+1})-\varphi(a^{n})}
\end{split}
\end{equation}
and hence
\begin{equation}\label{Eq:CLM}
\frac{M_{a^n, a^{n}}}{\varphi(a^{n+1})-\varphi(a^{n})}\to 0, \ \mathrm{as} \ n\to\infty.
\end{equation}

Now  by the assumption on $\varphi$,
\begin{equation}\label{phi-3}
\frac{\varphi(a^{n+1})}{\varphi(a^{n})}\leq C_2^{\lfloor \log_2 a\rfloor}
 \end{equation}
and by using Theorem 9.1 in Chobanyan, Leventhal and Mandrekar \cite{clm},
we see that (\ref{Eq:CLM}) implies
\begin{equation}
\lim_{n \to \infty} \frac{S_{0, n}} {\varphi(n)} = 0\qquad \hbox{a.s.}
\end{equation}
This finishes the proof of Theorem \ref{Th:SLLN}. \qed

The following is an immediate consequence of Theorem \ref{Th:SLLN} which is often
convenient to use.
\begin{corollary}
Let $\{\xi_{n},  n\ge 0\} $ be a sequence of random variables such that for
some $ 0<p<\infty$ and for all integers $k, n\in \N_0$
$$
\E |S_{k,n}|^p \leq g(n)
$$
for a numerical function $g$. If there is a nondecreasing function
$\varphi: \R_+\to \R_+$ as in Theorem \ref{Th:SLLN} such that
$$
\sum_{n=0}^{\infty}\frac{g(2^n)}{\varphi(2^{n})^p }<\infty,
$$
then $\{\xi_n, n \ge 0\}$ satisfies the $\varphi$-SLLN.
\end{corollary}

\section{Applications}

In this section, we show applications of Theorem \ref{Th:SLLN}
to quasi-stationary sequences of  random variables
and self-similar processes with stationary increments.

\subsection{Sequence of quasi-stationary random variables}

Strong laws of large numbers for sequences of quasi-stationary
random variables have been considered by several authors; see
\cite{Moricz77,clm} and the references therein. Let
$f:  \N_0 \to \R_+$ be a non-negative function. We say that a
real or complex-valued sequence $\{\xi_{n},\;n \in \N_0\}$ is
\textit{f-quasi-stationary} if $\E(\xi_k ) =0$,
$\E (\xi_{k}^2) < \infty$ for all
$ k \in \N_0$ and
$$\left|\E \big( \xi_{l}\bar{\xi}_{l+m}\big)\right| \leq f(m), \ \ \ \
l, m\;\in \N_0.$$

The following result refines Theorem 1 in \cite{Moricz77} and
Corollary 2 in \cite{clm}.
\begin{theorem}
Let $\{\xi_{n},\;n \in \N_0\}$ be an \textit{f-quasi-stationary} sequence
and let $\varphi$ be a non-decreasing function as in Theorem \ref{Eq:MomCon1}.
Define
$$
h(m)\equiv\sum_{n=\lfloor \log_{a}m\rfloor}^{\infty}\frac{a^n
}{{\varphi(a^n)^2}}.
$$
If $D :=\sum_{n=0}^{\infty}
\frac{a^n }{{\varphi(a^n)^{2}}}<\infty$ and
\begin{equation}\label{Eq:f-con}
\displaystyle
\sum_{m=1}^{\infty} f(m)h(m) < \infty,
\end{equation}
then the $\varphi$-SLLN holds for $\{\xi_{n}, \;n \in \N_0\}$.
\end{theorem}

\begin{proof} By the $f$-quasi-stationarity of $\{\xi_{n},\;n \in \N_0\}$,
we derive that for any $k, n \in \N_0$
$$
\E \bigg|\frac{S_{k, a^n}}{\varphi(a^n)}\bigg|^2 \leq \displaystyle
\sum_{m=0}^{a^n} \frac{f(m)(a^{n}-m)}{\varphi(a^n)^2} \leq
\frac{a^n}{{\varphi(a^n)^2}} \displaystyle \sum_{m=0}^{a^n} f(m).$$
It follows that
\begin{eqnarray*}
\displaystyle \sum_{n=0}^{\infty} \displaystyle \sup_{k \ge 0}\E
\bigg|\frac{S_{k, a^n}}{\varphi(a^n)}\bigg|^2
&\leq& \displaystyle \sum_{n=0}^{\infty} \sum_{m=0}^{a^n} \frac{a^n
f(m)}{{\varphi(a^n)^2}}\\
&\leq& D f(0)+ \displaystyle
\sum_{m=1}^{\infty} \sum_{n=\lfloor \log_{a}m\rfloor}^{\infty}
\frac{a^n f(m)}{{\varphi(a^n)^{2}}}\\
&\leq&D f(0)+ \displaystyle
\sum_{m=1}^{\infty}
 f(m)h(m)<\infty.
\end{eqnarray*}
Hence the conclusion follows from Theorem \ref{Eq:MomCon1}.
\end{proof}

\subsection{Self-similar processes with stationary increments}
Recall that a stochastic process $X= \{X(t), t \in \R_+\}$ with values in $\R$
is called a self-similar process with index $H > 0$ if for all constants $c > 0$,
\begin{equation}\label{Eq:SS}
\left\{X(ct), t \in \R_+\right\} \stackrel{d}{=} \, \left\{c^H X(t), t \in \R_+
\right\},
\end{equation}
where $\stackrel{d}{=}$ means equality of all finite dimensional distributions. $X$
is said to have stationary increments if for every $t_0\in \R_+$,
\begin{equation}\label{Eq:SI}
\{X(t_0 + t) - X(t_0), t \in \R_+\} \stackrel{d}{=} \, \{X(t) - X(0), t \in \R_+\}.
\end{equation}
If $X$ satisfies both (\ref{Eq:SS}) and \eqref{Eq:SI}, then we say that $X$ is $H$-SSSI.
We refer to Samorodnitsky and Taqqu \cite{ST94} and Embrechts and Maejima
\cite{EmbrechtsMaejima} for systematic accounts on self-similar processes.

The following theorem is concerned with asymptotic behavior of the sample function
$X(t)$ as $t \to \infty$.

\begin{theorem}\label{Coro:process}
Let $X= \{X(t), t \in \R_+\}$ be a real-valued $H$-SSSI process. If there is a constant
$p > \max\{1, 1/H\}$ such that $\E(|X(1)|^p) < \infty$, then for every $\ep > 0$,
\begin{equation}\label{Eq:coro1}
\lim_{t \to \infty} \frac{|X(t)|} {t^H (\log t)^{\frac 1 {p} + \ep}} = 0\qquad \hbox{a.s.}
\end{equation}
\end{theorem}

\begin{proof}\
For any integer  $n \ge 1$, define $\xi_n = X(n+1)-X(n)$. Then $\{\xi_n\}$ is stationary with
$\E(|\xi_n|^p) < \infty$. Let $\varphi(t) = t^H (\log t)^{\frac 1 {p} + \ep}$ and let $a \ge 2$
be an integer which satisfies the condition in Theorem 1.1.
By the $H$-SSSI property of $X$, we derive that
for all $k \ge 0$
\[
\E\big(|S_{k, a^n}|^p\big) = \E\big(|X(k + a^n) - X(k)|^p\big) = a^{nHp} \,\E(|X(1)|^p).
\]
This implies
\begin{equation}\label{Eq:SSSI-1}
\sum_{n=0}^\infty \sup_{k \ge 0} \E\bigg|\frac{S_{k, a^n}} {\varphi(a^n)}
\bigg|^p < \infty.
\end{equation}
Since $H > 0$,  we have $X(0) = 0$
a.s. \cite[p. 312]{ST94}. It follows from this fact and Theorem \ref{Th:SLLN}
that for every $\ep > 0$
\begin{equation}\label{Eq:coro2}
\lim_{n \to \infty} \frac{X(n)} {n^H (\log n)^{\frac 1 {p} + \ep}} = 0 \qquad \hbox{a.s.}
\end{equation}

To show (\ref{Eq:coro2}) still holds for continuous time $t$, we need an estimate on the tail
probability of $\max_{t \in [0, 1]}|X(t)|$. To this end, we note that for any integer $n\ge 1$
and $0 \le k \le n$,
\[
X\Big(\frac k n\Big) = \sum_{\ell=1}^k \bigg[X\Big(\frac \ell n\Big)- X\Big(\frac {\ell-1} n\Big)\bigg].
\]
Hence, for all $1 \le i < j \le n$, we have
\begin{equation}
\begin{split}
\E\Bigg(\bigg|\sum_{\ell=i}^j \bigg[X\Big(\frac \ell n\Big)-
X\Big(\frac {\ell-1} n\Big)\bigg]\bigg|^p\Bigg)
&= \E\bigg(\bigg|X\Big(\frac j n\Big) - X\Big(\frac {i-1} n\Big)\bigg|^p\bigg)\\
&= \bigg(\frac{j-i +1} n\bigg)^{Hp} \E\big(|X(1)|^p\big),
\end{split}
\end{equation}
where the second equality follows from the $H$-SSSI property of $X$. Since $p > \min\{1, 1/H\}$,
we see that the conditions of Theorem 3.1 in M\'oricz, Serfling and Stout \cite{MSS82} are
satisfied with
$$
g(i, j) = \big[\E\big(|X(1)|^p\big)\big]^{1/(Hp)}\,\frac {j-i +1} n
$$
which satisfies Condition (1.2) in \cite{MSS82}. It follows that
\begin{equation}\label{Eq:moment}
\E\bigg(\sup_{0 \le t \le 1}|X(t)|^p\bigg) \le K_1\, \E\big(|X(1)|^p\big) < \infty,
\end{equation}
where $K_1>0$ is an explicit constant depending on $H$ and $p$ only.
Combining (\ref{Eq:moment}) with the Markov inequality gives
\begin{equation}\label{Eq:MSS82}
\P\bigg(\sup_{0 \le t \le 1}|X(t)| \ge u\bigg) \le K_1\,\E\big(|X(1)|^p\big)\, u^{-p}.
\end{equation}
It follows from (\ref{Eq:MSS82}) and the $H$-SSSI property of $X$ that for any $\eta > 0$
\begin{equation}\label{Eq:Tail-SSSI}
\P\bigg(\sup_{a^n \le t \le a^{n+1}}|X(t) - X(a^n)| \ge \eta\, a^{nH} (\log a^n)^{1/p + \ep}\bigg)
\le K\, n^{-(1 + p \ep)}.
\end{equation}
Hence by the Borel-Cantelli lemma, we have
\begin{equation}\label{Eq:Tail-SSSI-2}
\lim_{n \to \infty} \frac{\sup_{a^n \le t \le a^{n+1}}|X(t) - X(a^n)|}
{a^{nH} (\log a^n)^{1/p + \ep} } = 0 \qquad \hbox{a.s.}
\end{equation}
It is clear that  (\ref{Eq:coro1}) follows from (\ref{Eq:coro2}) and (\ref{Eq:Tail-SSSI-2}).
\end{proof}

Applying a similar argument to self-similar stable processes with stationary increments, the
condition $H > 1/p$ can often be weakened, as shown by the next theorem.

\begin{theorem}\label{Coro:stable}
Let $X= \{X(t), t \in \R_+\}$ be a real-valued self-similar $\alpha$-stable process
with index  $H > 0$ and stationary increments. If the sample function $X(t)$ is almost
surely bounded on $[0, 1]$, then for every $\ep > 0$,
\begin{equation}\label{Th:coro-stable}
\lim_{t \to \infty} \frac{|X(t)|} {t^H (\log t)^{\frac 1 {\alpha} + \ep}} = 0\qquad \hbox{a.s.}
\end{equation}
\end{theorem}

\begin{proof}\
Let $a \ge 2$ be a fixed integer such that $a^{H \alpha} > \max\{2, 2^\a\}$. The same proof
as in that of Theorem \ref{Coro:process} shows that
\begin{equation}\label{Eq:lim-stable}
\lim_{n \to \infty} \frac{X(a^n)} {a^{nH} (\log a^n)^{1/\alpha + \ep}} = 0 \qquad {\rm a.s.}
\end{equation}
Under the assumption that $X$ has almost surely a bounded sample function on $[0, 1]$, we have
\[
\P\bigg(\sup_{0 \le t \le 1}|X(t)| \ge u\bigg) \le K\, u^{-\alpha}.
\]
See \cite[Theorem 10.5.1]{ST94}. This and the $H$-SSSI property of $X$ imply that
\begin{equation}\label{Eq:Tail-stable}
\P\bigg(\sup_{a^n \le t \le a^{n+1}}|X(t) - X(a^n)| \ge a^{nH} (\log a^n)^{1/\alpha + \ep}\bigg)
\le K\, n^{-(1 + \a \ep)}.
\end{equation}
Hence (\ref{Th:coro-stable}) follows from (\ref{Eq:lim-stable}), (\ref{Eq:Tail-stable})
and the Borel-Cantelli lemma.
\end{proof}

As an example, let us consider the linear fractional stable motion. Given constants
$\alpha \in (0, 2)$ and $H \in (0, 1)$, the $\alpha$-stable process
$\{L_{\alpha,H}(t), t \in \R\}$ defined by
\begin{equation}\label{Def:LFSS}
L_{\alpha,H}(t)=\int_{\R}\left[(t-s)_+^{H-1/\alpha}-(-s)_+^{H-1/\alpha}\right]\,A(ds)
\end{equation}
is called a linear fractional stable motion (LFSM) with indices $\alpha$ and $H$.
In the above, $a_+ = \max\{0, a\}$ for all $ a \in \R$ and $\{A(t), t \in \R\}$ is
a two-sided strictly stable L\'evy process
of index $\alpha$. Note that, when $H =  1/\alpha$, $L_{\alpha,H}(t) = A(t)$ for
all $t \ge 0$. When  $H \ne 1/\alpha$, the stochastic integral in (\ref{Def:LFSS})
is well-defined because
\[
\int_{\R}\Big|(t-s)_+^{H-1/\alpha}-(-s)_+^{H-1/\alpha}\Big|^\alpha\,dr
< \infty.
\]
See \cite[Chapter 3]{ST94}. By (\ref{Def:LFSS}), it can be verified that
$\{L_{\alpha,H}(t), t \in \R\}$ is
$H$-SSSI \cite[Proposition 7.4.2]{ST94}. It
is an $\alpha$-stable analogue of fractional Brownian motion and its
probabilistic and statistical properties have been investigated by many authors.
In particular, it is known that
\begin{itemize}
\item[(i)]\ If $1/\alpha<H<1$ (this is possible only when $1 < \alpha < 2$),
then the sample
function of $\{L_{\alpha,H}(t), t \in \R\}$ is almost surely continuous.
\item[(ii)]\ If $ 0< H<1/\alpha$, then the the sample function of $\{L_{\alpha,H}(t),
t \in \R\}$ is almost surely unbounded on every interval of positive length.
\end{itemize}
We refer to \cite[Chapters 10 and 12]{ST94} and
Takashima \cite{Taka89} for more information on asymptotic properties of LFSM.

The following is a corollary of Theorem \ref{Coro:stable}.
\begin{itemize}
\item
If $1/\alpha \le H <1$, then for every $\ep > 0$,
\begin{equation}\label{Eq:LFSS}
\lim_{t \to \infty} \frac{|L_{\alpha, H}(t)|} {t^H (\log t)^{\frac 1 {\a} + \ep}} = 0
\qquad \hbox{a.s.}
\end{equation}
\item
If $0 <H<1/\alpha$, then for every $\ep > 0$,
\begin{equation}\label{Eq:LFSS2}
\lim_{n \to \infty} \frac{|L_{\alpha, H}(n)|} {n^H (\log n)^{\frac 1 {\a} + \ep}} = 0
\qquad \hbox{a.s.}
\end{equation}
\end{itemize}

For $1/\alpha < H<1$, (\ref{Eq:LFSS}) is proved in Ayache, Roueff and
Xiao \cite{ARX07} by using the wavelet methods. For $0 <H<1/\alpha$, even though
$\{L_{\alpha,H}(t), t \in \R\}$ is almost surely unbounded on every interval
of positive length, its limiting behavior along a fixed sequence is still
similar to the  $1/\alpha \le H<1$ case.








\end{document}